\documentclass[final,3p,twocolumn]{elsarticle}
\pdfoutput=1
\usepackage{graphicx,subfigure,url}
\usepackage{amsmath,amsfonts,color}
\usepackage[pdftex,colorlinks]{hyperref}
\definecolor{darkblue}{cmyk}{1,0,0,0.8}
\definecolor{darkred}{cmyk}{0,1,0,0.7}
\hypersetup{anchorcolor=black,
  citecolor=darkblue, filecolor=darkblue,
  menucolor=darkblue,pagecolor=darkblue,urlcolor=darkblue,linkcolor=darkblue}

\usepackage{typearea}
\typearea{23}
\usepackage[charter]{mathdesign}
\usepackage{microtype}

\renewcommand{\epsilon}{\varepsilon}
\renewcommand{\epsilon}{e}
\newlength{\picrefwidth}
\journal{International Journal of Nonlinear Mechanics}
\begin{document}
\begin{frontmatter}
  \title{Dynamics of the Nearly Parametric Pendulum}
  \author[AB]{B. Horton}
  \author[AB]{J. Sieber}
  \author[AB]{J. M. T. Thompson}
  \author[AB]{M. Wiercigroch\corref{cor}}
  \address[AB]{Centre for Applied Dynamics Research\\
    University of Aberdeen, School of Engineering,
    Kings College, Aberdeen,\\ AB24 3UE, UK}
  \cortext[cor]{Corresponding author.}
  \ead{m.wiercigroch@abdn.ac.uk}
  \begin{abstract}
    Dynamically stable periodic rotations of a driven pendulum provide
    a unique mechanism for generating a uniform rotation from bounded
    excitations.  This paper studies the effects of a small
    ellipticity of the driving, perturbing the classical parametric
    pendulum.  The first finding is that the region in the parameter
    plane of amplitude and frequency of excitation where rotations are
    possible increases with the ellipticity.  Second, the resonance
    tongues, which are the most characteristic feature of the
    classical bifurcation scenario of a parametrically driven
    pendulum, merge into a single region of instability.
  \end{abstract}
\end{frontmatter}
\section{Introduction} \label{sec:intro}

The driven pendulum is a generic model used for studying nonlinear
dynamics in mechanics \cite{RefWorks:60} and beyond
\cite{RefWorks:83,mouchet1,RefWorks:69,AJP000525,RefWorks:66}. Its
geometric nonlinearity can be modeled reliably (in contrast to other
nonlinear effects such as friction), and mechanical pendula are
amenable to experimental investigations. The dynamical properties of
the classical parametrically driven pendulum, such as resonances,
escape from a potential well, symmetry-breaking, and periodic and
chaotic attractors, have been explored in detail experimentally
\cite{RefWorks:60,AJP000909,AJP000821,magneticpendulum,BHpaper1,xuzamm}
and theoretically
\cite{RefWorks:51,RefWorks:63,RefWorks:65,xuwierci,lenci2,Isohatala,%
RefWorks:67,RefWorks:52,stewartharmonicbalance,RefWorks:49,RefWorks:62,%
smithblackburncircuit,RefWorks:53,Kobes1,Kim1,RefWorks:72}.

This paper studies what happens to the well-studied bifurcation
scenarios of the parametrically excited pendulum if the driving of the
pivot of the pendulum follows a narrow upright ellipse; see
figure~\ref{fig:parapendulum}.  One motivation for studying elliptic
excitation is that only the elliptic component of an arbitrarily
shaped periodic excitation has an effect on a rotating pendulum for
large excitation frequencies; see section~\ref{sec:reduce} for an
explanation.  Moreover, elliptic excitation is typical if the pendulum
base is floating on water waves: a small freely floating body moves
along an ellipse. This effect is similar to the elliptic motion of an
off-center surface point of a plate excited by a circular traveling
bending wave (a principle that is exploited in rotary ultrasonic
motors \cite{HW92, HM95}).

We say that the pendulum \emph{rotates} if the long-time average of
the angular velocity is non-zero.  Stable periodic rotations occur
naturally over a large range of excitation parameters in the
parametrically driven pendulum \cite{RefWorks:51}.  Thus, a rotating
pendulum provides a unique mechanism for generating a uniformly
one-directional rotation from a bounded excitation. This is a
potential physical principle for harnessing the energy of vibrations,
which are not necessarily purely in the vertical direction.  The other
motivation for focusing on rotating attractors is that the rotating
pendulum is ideal for developing and testing non-invasive bifurcation
and chaos control methods \cite{OGY90,P92} in a real experiment:
periodic rotations are reliably controllable by superimposing feedback
control onto the excitation without changing the shape of the
excitation. This is not true in general for small-amplitude
oscillations around the hanging-down position \cite{WW00}.

The dimensionless equation of motion for the elliptically excited
pendulum is
\begin{equation}
  \label{eq:elliptic1a}
  \ddot{\theta}+\gamma\dot{\theta}+\left(1+p\cos\left(\omega
      t\right)\right)\sin\theta+e p\sin(\omega t)\cos\theta=0
\end{equation}
where $\gamma$ is the dimensionless viscous damping coefficient, $p$
is the scaled excitation amplitude, $\omega$ the rescaled excitation
frequency, and $e$ is the ratio between the horizontal and the
vertical diameter of the upright ellipse traced out by the pivot
during each period (see figure~\ref{fig:parapendulum}). The classical
parametrically excited pendulum corresponds to the setting $e=0$.

The two main effects of a small non-zero ellipticity $e$ of the
excitation are:
\begin{enumerate}
\item\label{tongues} The classical resonance tongues for the 1:2 and
  the 1:1 resonance of the parametrically excited pendulum \cite{RefWorks:51}
  merge into a single region of instability of the small-amplitude
  period-one libration around the hanging-down position of the pendulum.
\item\label{pref} If the ellipticity $e$ is non-zero the pendulum is
  no longer symmetric with respect to reflection
  $\theta\mapsto-\theta$, which causes a preference for rotations that
  have the same direction as the motion of the pivot.  Effectively,
  the range of possible excitation frequencies and amplitudes where
  rotations are supported increases for increasing ellipticity. The
  preferred direction of rotation has the same sense (clockwise or
  anti-clockwise) as the motion of the base around the ellipse (for
  example, clockwise rotation is preferred if the pivot moves
  clockwise around the ellipse) because this rotation picks up energy
  from the additional excitation in the horizontal direction.
\end{enumerate}
These two observations are, in short, the key findings of the paper.
Point~\ref{pref} is, for large frequencies, universal for all shapes
of excitation that have a dominant vertical component.  This will be
shown in section~\ref{sec:reduce} by averaging the equation of motion
for a pendulum with arbitrary periodic
excitation. Section~\ref{sec:reduce} also gives an approximate
expression for the onset of rotations that is valid for all shapes of
excitation if the forcing frequency $\omega$ is
large. Section~\ref{sec:model} shows how the non-dimensionalized
equation of motion \eqref{eq:elliptic1a} is related to the original
equation of motion describing a physical pendulum driven by a slider
along an ellipse. Section~\ref{sec:bifoverview} and
section~\ref{sec:rot} give two-parameter overviews of changes to
the classical structure of resonance tongues and to the existence
regions of rotations.  Section~\ref{sec:1d} shows one-parameter
bifurcation diagrams (for increasing forcing amplitude) to illustrate
how the different attractors are connected by the bifurcations shown
in the figures \ref{fig:paraspaces_oscs} and
\ref{fig:paraspaces_rots}.

\section {Rotations in the high-frequency forcing regime}
\label{sec:reduce}
Let us assume that the pivot of the pendulum is driven periodically
with high frequency along an arbitrary path. Then the inclination
angle $\theta$ of the pendulum is governed by the equation of motion
\begin{equation}
  \label{eq:gen}
  \ddot\theta+\gamma\dot\theta+\omega f_y(t)\sin\theta+\omega 
  f_x(t)\cos\theta=0\mbox{,}
\end{equation}
where the two coordinates of the force on the pendulum bob caused by
the displacement of the pivot and by gravity, $\omega f_x(t)$ and
$\omega f_y(t)$, have period $2\pi/\omega$. The high-frequency regime
is the parameter range where the frequency $\omega$ is large. The
effect of high-frequency elliptic excitation on small amplitude
oscillations has been studied analytically recently using averaging
techniques \cite{FT08}. In the high-frequency regime the excitation
forces are typically large, even if the excitation amplitude of the
pivot (for example, the size of the ellipse in
figure~\ref{fig:parapendulum}) is small. Thus, we have put the scaling
factor $\omega$ expressly in front of $f_x$ and $f_y$, and assume that
$\gamma$, $|f_x(t)|$ and $|f_y(t)|$ are of order $1$ (this corresponds
to an excitation amplitude of order $\omega^{-1}$ for the
pivot). Rotations of the pendulum in the positive direction
($\dot\theta>0$) with a frequency close to the forcing frequency are
solutions of \eqref{eq:gen} for which the quantity
$\phi(t)=\theta(t)-\omega t$ is bounded for all times. We insert
$\phi$ into equation \eqref{eq:gen} and average \eqref{eq:gen} to
second order over one period. The second-order averaged equation for
$\phi$ is valid on the slow time scale $\sqrt{\smash[b]{\omega}}t$
(which is faster than the original time scale $t$ of \eqref{eq:gen}
but slower than the time scale of the forcing $\omega t$):
\begin{equation}\label{eq:genavg}
  \ddot \phi+\frac{\gamma}{\sqrt{\smash[b]{\omega}}}\,\dot\phi+\gamma+
  \frac{1}{2}\left[f_y^c-f_x^s\right]\sin\phi+
  \frac{1}{2}\left[f_y^s+f_x^c\right]\cos\phi=0
\end{equation}
where
\begin{align}
  \label{eq:fcoeffs}
  f_{x,y}^c&=\ \frac{\omega}{\pi}\int_0^{2\pi/\omega}
  f_{x,y}(s)\cos(\omega s)\,\mathrm{d}s\mbox{,}\nonumber\\
  f_{x,y}^s&=\frac{\omega}{\pi}\int_0^{2\pi/\omega}
  f_{x,y}(s)\sin(\omega s)\,\mathrm{d}s
\end{align}
are the coefficients of the first Fourier modes of $f_x$ and
$f_y$. All coefficients in \eqref{eq:genavg} are at most of order $1$,
and the periodic terms that have been dropped in the averaging
procedure are of order $\omega^{-1}$. Thus, for a large frequency
$\omega$, only the first Fourier coefficients of the excitation,
$f_{x,y}^{c,s}$, have an influence at the leading order.  We can
assume that one of the four leading Fourier coefficients is zero
without loss of generality (we can shift time to make, for example,
$f_y^s=0$), and introduce three parameters to describe the other three
coefficients:
\begin{equation}
  \label{eq:ffourierc}
    \omega f_y^c=p\mbox{,}\quad
    f_y^s=0\mbox{,}\quad
    \omega f_x^c=\epsilon p\cos\alpha\mbox{,}\quad
    \omega f_x^s=\epsilon p\sin\alpha\mbox{.}
\end{equation}
where $p>0$ can be large (of order $\omega$) and
$\alpha\in[0,2\pi]$. The parameter $|\epsilon|$ is the ratio between
$\sqrt{\smash[b]{(f_x^c)^2+(f_x^s)^2}}$ and $f_y^c$, and the parameter
$\alpha$ describes the phase shift between the horizontal and the
vertical component of the first harmonic of the forcing. Using these
parameters the averaged equation \eqref{eq:genavg} for $\phi$ becomes
\begin{equation}\label{eq:fastrot}
  \ddot \phi+\frac{\gamma}{\sqrt{\smash[b]{\omega}}}\,\dot\phi+\gamma+
  \frac{p}{2\omega}\left[\left(1-\epsilon\sin\alpha\right)\sin\phi+
    \epsilon\cos\alpha\cos\phi\right]=0\mbox{.}
\end{equation}
Positively directed rotations of period $2\pi/\omega$ correspond
approximately to equilibria of \eqref{eq:fastrot} in the following
sense: if $\gamma$, $e$ and $p/\omega$ are at most of order one, and
the averaged equation \eqref{eq:fastrot} has an equilibrium $\phi_0$
then the original forced equation \eqref{eq:gen} has a solution
$\theta$ satisfying for all times $t$
\begin{equation}
  \label{eq:avgjust}
  \theta(t)-\omega t=\phi_0+r(t)
\end{equation}
where $|r(t)|\ll1$ and $r$ has period $2\pi/\omega$. The stability
properties of the equilibrium $\phi_0$ also transfer to the rotation
$\theta$: if $\phi_0$ is stable then $\theta$ is stable, if $\phi_0$
is a saddle then $\theta$ is a rotation of saddle-type. Moreover,
bifurcations of the equilibria of \eqref{eq:fastrot} are also
transferred: since \eqref{eq:fastrot} is dissipative, only saddle-node
bifurcations can occur. Indeed, if $2\gamma\omega= p\sqrt{\smash[b]{1+
    \epsilon^2-2\epsilon\sin\alpha}}$ then the averaged equation
\eqref{eq:fastrot} has a saddle-node bifurcation, which implies that
the original system \eqref{eq:gen} has a saddle-node bifurcation of
rotations at parameters nearby.  If we replace $\alpha$ by $-\alpha$
in \eqref{eq:fastrot} then the equilibria of \eqref{eq:fastrot}
correspond to periodic rotations in the negative direction (that is,
to solutions $\theta$ of \eqref{eq:gen} satisfying
$\theta(t+2\pi/\omega)=\theta(t)-2\pi$ for all $t$), and the sign in
front of $2\epsilon\sin\alpha$ changes to $+$ in the condition for the
saddle-node bifurcation.  Thus, for large frequency $\omega$, periodic
rotations of \eqref{eq:gen} satisfying
$\theta(t+2\pi/\omega)=\theta(t)\pm2\pi$ are born in a saddle-node
bifurcation defined (up to terms of order $\omega^{-1}$ by
\begin{equation}
  \label{eq:foldapp}
  2\gamma= \frac{p}{\omega}\sqrt{\smash[b]{1+
      \epsilon^2\mp2\epsilon\sin\alpha}}
\end{equation}
if $\gamma$ and $\epsilon$ are of order $1$ and $p$ is of order
$\omega$.  One of the rotations emerging from the saddle-node is
stable and remains stable for arbitrarily large $p$ as long as the
averaging approximation is valid. This implies that the curve of
period doublings of rotations which we observed numerically (see
figure~\ref{fig:paraspaces_rots} in section~\ref{sec:results}) has to
grow super-linearly in $p$ for increasing $\omega$.  In summary, for
large frequency $\omega$, we have:
\begin{enumerate}
\item If the force amplitude is at most of order $\omega$ then the
  existence and stability of rotations is entirely determined by the
  first Fourier mode of the excitation.
\item A stable and a saddle-type periodic rotation satisfying
  $\theta(t+2\pi/\omega)=\theta(t)\pm 2\pi$ are born in a saddle-node
  bifurcation near parameter values given by \eqref{eq:foldapp} where
  $p$, $e$ and $\alpha$ define the amplitudes of the first Fourier
  mode as given in \eqref{eq:ffourierc}.
\item The stable rotation remains stable for increasing $p$ over a
  large region of parameter values of $p$ (as long as $p$ is of
  similar magnitude to $\omega$ and the averaged equation is a valid
  approximation).
\item The difference between the positive and negative directions of
  rotation is maximal for $\alpha=\pm\pi/2$ in expression
  \eqref{eq:foldapp} for the onset of rotations. This corresponds to
  the case where the first harmonics form an upright ellipse.
\end{enumerate}
We note that one can extend the averaging technique to frequencies
$\omega$ of order $1$: for small damping $\gamma$ and small forcing
$p$ (and $\epsilon p$) one can average along the integral curves of
the unforced and undamped pendulum. This technique was used in
\cite{lenci2} for the model of a parametrically driven pendulum
($e=0$) and can be applied also for a forcing of general harmonic
shape (such as the elliptically driven pendulum). Using this refined
averaging we found that for small damping ($\gamma=0.1$) the
expression \eqref{eq:foldapp} for the saddle-node bifurcation of
rotations is a good approximation for $\omega>1.5$ (if all quantities
refer to the non-dimensionalized equation \eqref{eq:elliptic1a} where
$\alpha=\pi/2$).

\section{Modelling of the elliptically driven pendulum}
\label{sec:model}
The numerical results in section~\ref{sec:results} discuss what
happens if the forcing of the pendulum deviates from the classical
parametrically driven pendulum and the forcing frequency is near the
main resonance tongues known from the parametric case. We excite the
pendulum harmonically along a narrow ellipse.  This corresponds to a
choice of
\begin{equation}
  \label{eq:ffourier}
    \omega f_{x}(t)=\epsilon p\cos(\omega t-\alpha)\mbox{,}\quad
    \omega f_{y}(t)=p\cos(\omega t)+1
\end{equation}
in \eqref{eq:gen}. In \eqref{eq:elliptic1a} and \eqref{eq:ffourier} we
use the convention that the angle $\theta=0$ corresponds to the
hanging-down position of the pendulum such that the force due to
gravity contributes a positive constant term to the coefficient
$\omega f_y$ in front of $\sin\theta$ but nothing to $\omega f_x$.  If
we assume that the vertical component of the forcing is dominant then
$|\epsilon|$ is significantly less than $1$ such that the overall
forcing amplitude $p\sqrt{\smash[bt]{1+\epsilon^2}}$ is controlled to
first order of $\epsilon$ by $p$ only since
\begin{displaymath}
  \frac{\mathrm{d}}{\mathrm{d}\epsilon}
  \left[\sqrt{\smash[b]{1+\epsilon^2}}\right]_{e=0}=0\mbox{,}
\end{displaymath}
and $\epsilon$ controls the ellipticity.

The parameter $\alpha$ controls the inclination of the ellipses
ranging between excitation along a straight line ($\alpha=0$) and the
family of upright ellipses ($\alpha=\pm\pi/2$). Note, however, that
$\alpha$ is not identical to the inclination angle of the ellipse: for
example, for $\alpha=0$ the inclination of the straight line is
determined by $e$.
\begin{figure}[ht]
    \begin{center}
    {\includegraphics[scale=1]{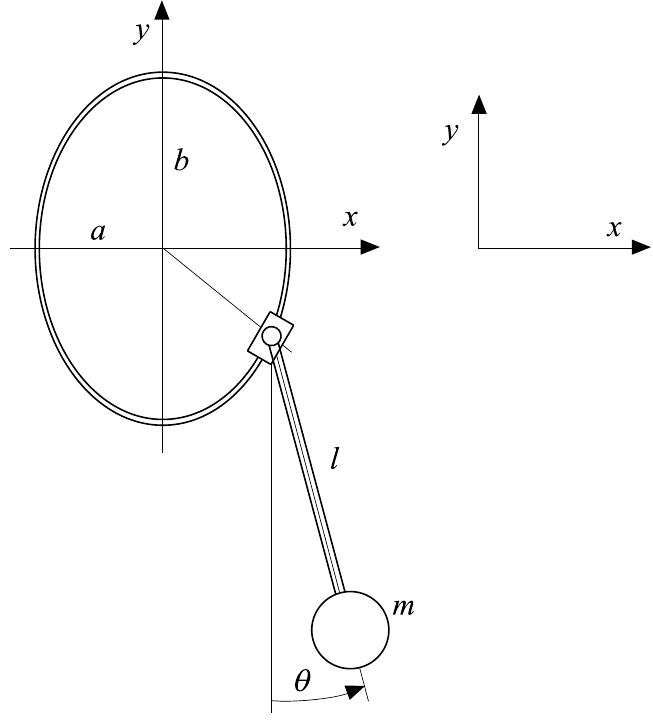}}
    \end{center}
    \caption{Schematic of an elliptically excited mechanical pendulum.}
    \label{fig:parapendulum}
\end{figure}

The extreme case of parametric excitation corresponds to
$(\epsilon,\alpha)=(0,0)$.  The other extreme case of horizontal
excitation, which would be singular ($\epsilon\to\infty$) with our
choice of parameters, has been studied theoretically
in~\cite{horizontalpend,horizontalpend2}. 

The approximate expression \eqref{eq:foldapp} for the onset of
rotations shows that the difference betwen both directions of rotation
is most prominent if the parametric excitation is perturbed into an
upright ellipse ($\alpha=\pi/2$ in \eqref{eq:ffourier}). Thus, we
restrict our numerical study in section~\ref{sec:results} to the
specific model \eqref{eq:elliptic1a}, which corresponds to
$\alpha=\pi/2$. This gives rise to equation~\eqref{eq:elliptic1a} for
the elliptically excited pendulum as proposed in the introduction.

Figure~\ref{fig:parapendulum} shows a mechanical representation of
this model: a pendulum having mass $m$ and length $l$ is driven by a
slider on an elliptic kinematic constraint. The slider is connected to
the pendulum rod via a pin joint.  The non-dimensional parameters and
the non-dimensional time of model \eqref{eq:elliptic1a} can be
obtained from the corresponding quantities of the mechanical
representation by the scaling
\begin{equation}\label{dimensionalization}
\begin{aligned}
 \gamma&=\frac{c}{\omega_0ml^2}\mbox{,} &
 \omega&=\frac{\Omega}{\omega_0}\mbox{,}&
 t_\mathrm{scaled}&=\omega_0\, t_\mathrm{physical}\mbox{,}\\
\omega_0&=\sqrt{\frac{g}{l}}\mbox{,}&
 p&=\frac{a\Omega^2}{g}\mbox{,}&
 ep&=\frac{b\Omega^2}{g}\mbox{.}
\end{aligned}
\end{equation}
In \eqref{dimensionalization} $a$ is the amplitude of the vertical
displacement excitation $a\cos(\Omega t)$, $b$ is the amplitude of
the horizontal displacement excitation $b\sin(\Omega t)$, $\Omega$ is the
driving frequency, $g$ is the acceleration due to gravity,
$\omega_0$ is the linear natural frequency of the pendulum at the
hanging-down angle $\theta=0$, $l$ is the length of the (mass-less)
pendulor arm, $m$ is the mass of the pendulum bob, and $c$ is the
viscous damping coefficient in the mechanical representation shown in
figure~\ref{fig:parapendulum}.

\section{Resonance Structure}\label{sec:results}

In this section we analyse how the introduction of a nonzero
ellipticity $\epsilon$ changes the resonance structure by constructing
two-parameter bifurcation diagrams in the
$\left(\omega,p\right)$-plane. We also use cross-sections of these
diagrams (one-parameter bifurcation diagrams) at constant frequencies
$\omega$ to make the connection between the different dynamical
regimes visible. Throughout our study the bifurcation parameters are
the excitation frequency $\omega$, the excitation amplitude $p$ and
the ellipticity $e$ of the excitation.  The parameter $e$ perturbs the
reflection symmetry of the parametrically driven pendulum in the
following way: if $(\theta(t),\dot\theta(t))$ is a solution for $e$
then $(-\theta(t),-\dot\theta(t))$ is a solution for $-e$. Thus, for
any ellipticity $e$ the bifurcations obtained for $e$ are identical to
those obtained for $-e$. This implies that we can restrict our
attention to $e\geq0$.

We slice the three-dimensional $(\omega,p,e)$-space along three
two-dimensional planes by constructing three two-parameter bifurcation
diagrams in the $(\omega,p)$-plane for three different values of $e$:
$e=0$, $e=0.1$ and $e=0.5$. The case $e=0$ is the classical
parametrically driven pendulum as studied in
\cite{RefWorks:51,RefWorks:49,xuwierci}.  The case $e=0.1$
shows how a small perturbation of the reflection symmetry affects the
classical bifurcation scenario and $e=0.5$ provides a picture of how
the bifurcation scenario changes as the system deviates further from
the parametric pendulum case towards a circular excitation.

The only remaining parameter in the non-dimensionalized equation
\eqref{eq:elliptic1a} is the dimensionless viscous damping
$\gamma$. We choose
\begin{math}
  \gamma=0.1
\end{math}
to make our results comparable with the results of the
previous studies \cite{RefWorks:51,RefWorks:49,xuwierci}.

\begin{figure}[ht]
  \begin{center}
    {\includegraphics[width=1\picrefwidth]{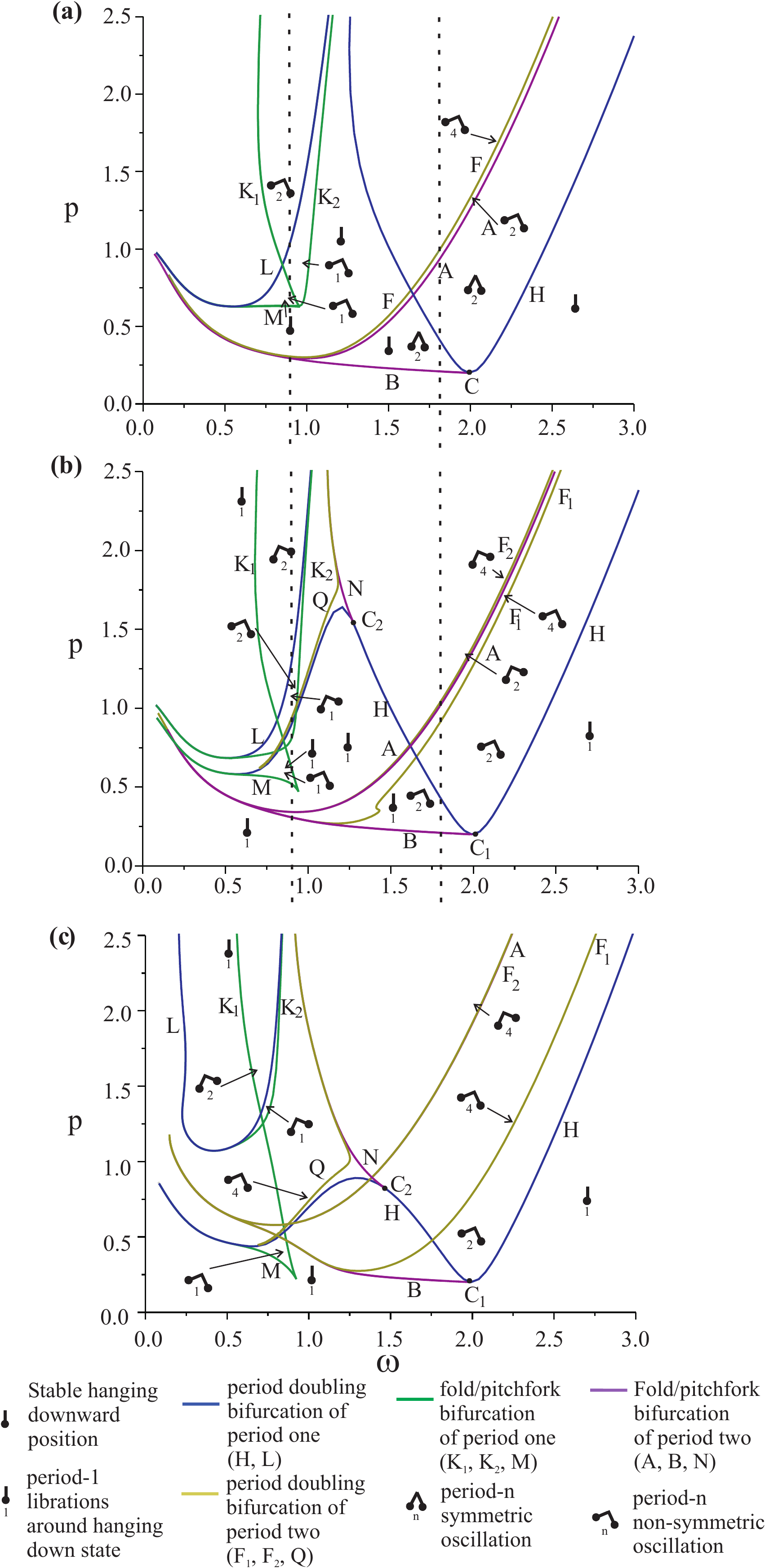}}
  \end{center}
  \caption{(colour online) Bifurcation diagrams in the
    $(\omega,p)$-plane for (a) $e=0$, (b) $e=0.1$, (c) $e=0.5$,
    computed with (R)AUTO \cite{auto,S07}. The figure
    \ref{fig:oneparabif1} shows the one-parameter bifurcation diagrams
    along the dashed lines in (a) and (b) for $\omega=0.87$ and
    $\omega=1.8$, respectively. The damping $\gamma$ is $0.1$.}
  \label{fig:paraspaces_oscs}
\end{figure}
In the following we will discuss the bifurcations of oscillations and
rotations separately.  Oscillations are periodic orbits that stay in
the potential well of the undriven conservative pendulum around the
hanging-down position $\theta=0$. The average of the angular velocity
over one period of an oscillation is zero. Rotations leave this
potential well and have a non-zero average angular velocity along one
period (for period-one rotations the average angular velocity is $\pm
\omega$). 
We will present
rotations and oscillations always in separate figures because they
coexist over large parameter ranges and there is no local bifurcation
linking the two types of periodic orbits.
\begin{figure}[ht]
 \begin{center}
         {\includegraphics[width=1\picrefwidth]{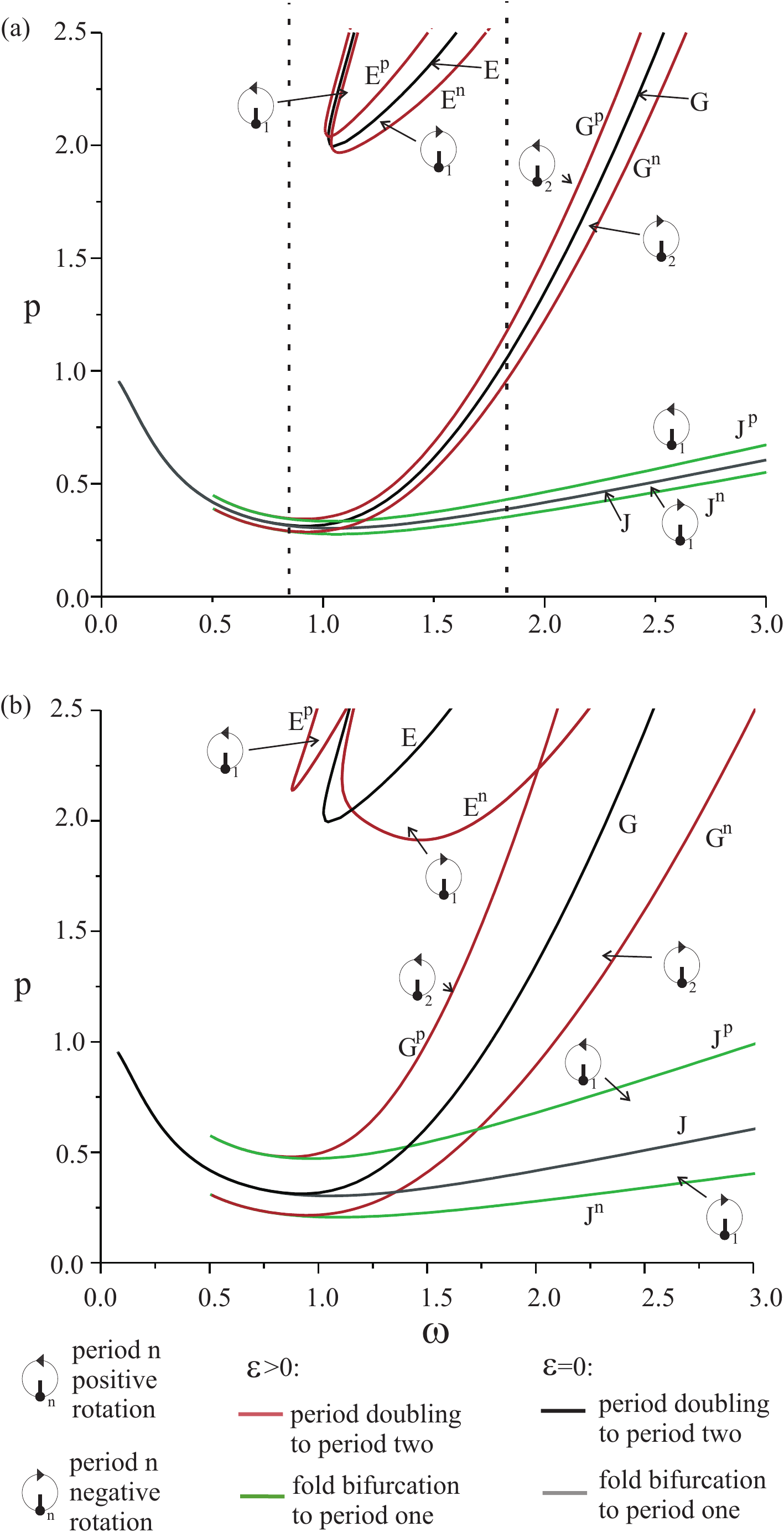}}
     \end{center}
     \caption{(colour online) Bifurcation diagrams in the
       $(\omega,p)$-plane comparing the bifurcations of the period-one
       rotations for $e=0.1$ (panel(a)) and $e=0.5$ (panel (b)) with
       the classical case $e=0$, computed with (R)AUTO
       \cite{auto,S07}.  The figure \ref{fig:oneparabif1} shows the
       one-parameter bifurcation diagrams along the dashed lines in
       (a) for $\omega=0.87$ and $\omega=1.8$, respectively.  The
       damping $\gamma$ is $0.1$.}
     \label{fig:paraspaces_rots}
\end{figure}

\subsection{Overview of oscillations in the ($\omega$,p) plane}
\label{sec:bifoverview}
Figure \ref{fig:paraspaces_oscs} shows the two-parameter bifurcation
diagrams in the $(\omega,p)$-plane for oscillations.  Panel (a) shows
the classical diagram for $e=0$, panel (b) presents the diagram for
$\epsilon=0.1$, and panel (c) presents the diagram for $\epsilon=0.5$.
The symbols between bifurcation curves in
figure~\ref{fig:paraspaces_oscs} indicate which attractors are
observable in the different regions. The most prominent features of the
classical diagram \ref{fig:paraspaces_oscs}(a) are the two main
resonance tongues where the hanging-down position $\theta=0$ loses its
stability: the 1:2 resonance at $\omega\approx 2$ and the 1:1
resonance at $\omega\approx 1$. The 1:2 resonance tongue is bounded by
the blue period-doubling curve H, and the 1:1 resonance tongue, which
starts at a larger value of forcing ($p\approx0.7$), is bounded by a
pitchfork bifurcation curve (light green curves K$_{1,2}$ in figure
\ref{fig:paraspaces_oscs}(a)).  Both tongues are separated by a region
of stability of $\theta=0$ between the curves K$_2$ and H in
figure~\ref{fig:paraspaces_oscs}(a).  The period doubling H bounding
the 1:2 resonance has a degeneracy at the point C: it is supercritical
to the right of C and subcritical to the left of C.

The most significant change for nonzero ellipticity $e$ is that the
two resonance tongues merge into a single region of instability. The
period-doubling curve H merges with one of the non-symmetric period
doubling curves L of the 1:1 tongue. This period doubling (still
called H in figures \ref{fig:paraspaces_oscs}(b) and
\ref{fig:paraspaces_oscs}(c)) and the fold curve K$_1$ form the
stability boundary for the small-amplitude libration of period one
around $\theta=0$, which is a perturbation of order $\epsilon$ of the
hanging-down equilibrium position $\theta=0$ of the classical
parametrically driven pendulum ($\epsilon=0$). The period doubling is
subcritical between the points C$_1$ and C$_2$ along the curve H.

The one-parameter bifurcation diagrams along the parameter paths
marked as dashed lines in figure~\ref{fig:paraspaces_oscs}(a) and (b)
are discussed in detail in section~\ref{sec:1d}. They show how
the other bifurcation curves in
figure~\ref{fig:paraspaces_oscs} form the stability boundaries for the
more complex oscillations. The values $\omega=1.8$ and $\omega=0.87$
for these parameter paths are representative for the 1:2 and the 1:1
resonance, respectively. They are the same as in \cite{RefWorks:49},
which studied the parametric case $e=0$.

\subsection{Overview of rotations in the ($\omega$,p) plane}
\label{sec:rot}
Figure~\ref{fig:paraspaces_rots} shows the bifurcations of period-one
rotations for $e=0.1$ (figure~\ref{fig:paraspaces_rots}(a)) and
$e=0.5$ (figure~\ref{fig:paraspaces_rots}(b)). The bifurcations of
rotations in the parametric case $e=0$ are included (J and G, in black
and grey) in both panels to show the effect of the nonzero ellipticity
$e$. In all cases the stable period-one rotations are born, for
increasing forcing $p$, in a fold bifurcation (curves J, J$^p$ and
J$^n$ in figure~\ref{fig:paraspaces_rots}) and lose their stability in
a period doubling bifurcation (curves G, G$^p$ and G$^n$ in
figure~\ref{fig:paraspaces_rots}) as $p$ increases further. For even
higher forcing the period-one rotation regains its stability (in the
period doublings E, E$^p$ and E$^n$).

The most notable effect of the nonzero ellipticity is that all
bifurcations are shifted toward lower forcing for rotations in the
negative direction (that is, in the same rotiational sense as the
motion of the base along the ellipse). The bifurcations of rotations
in the positive direction are shifted upward. Approximation
\eqref{eq:foldapp} estimates this effect in the limit of high
frequency.

According to figure~\ref{fig:parapendulum} the rotation in the negative
direction ($\dot\theta<0$) rotates in the same direction as the base
of the pendulum, corresponding to $\alpha=-\pi/2$ in
\eqref{eq:foldapp}.  Thus, for $e>0$ the curve J$^n$ is shifted
downwards from J by $2\gamma\omega e/(1+e)$ and the curve J$^p$ is
shifted upwards from J by $2\gamma\omega e/(1-e)$ in the limit
$\omega\to\infty$. The equilibria of the averaged equation
\eqref{eq:fastrot} show that negative rotations ($\alpha=-\pi/2$) pick
up energy from the horizontal component of the forcing on average
(that is, the factor in front of $\sin\phi$ is larger than one)
whereas the positive rotations lose energy.

\subsection{One-parameter diagrams for varying forcing amplitude}
\label{sec:1d} 
\begin{figure}[ht]
 \begin{center}
   {\includegraphics[width=1\picrefwidth]{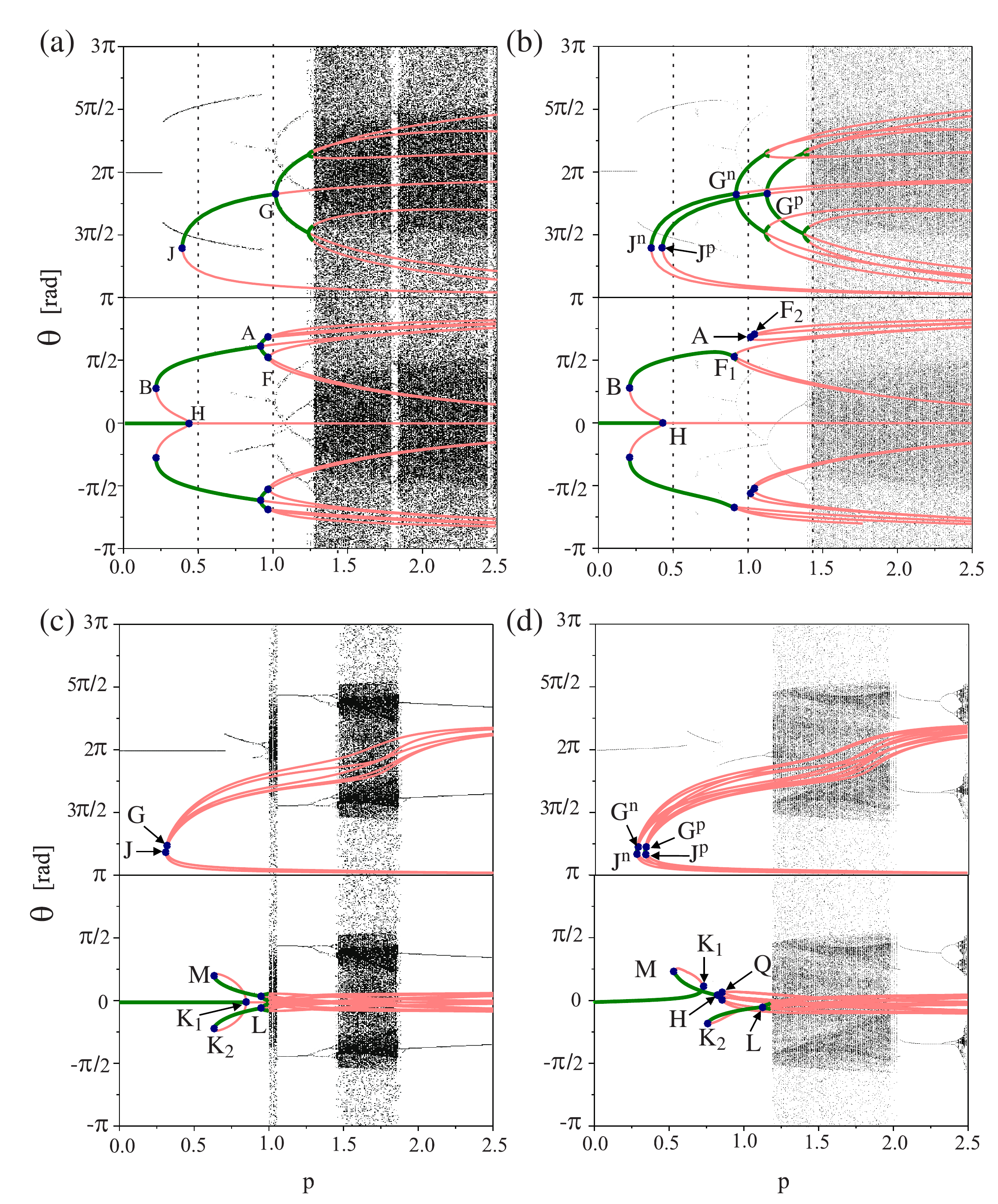}}
 \end{center}
 \caption{(colour online) One parameter bifurcation diagrams for the
   stroboscopic map of \eqref{eq:elliptic1a} for different values of
   $\omega$, varying the scaled forcing amplitude $p$ as the
   bifurcation parameter. The other parameters are $e=0$, $\omega=1.8$
   in panel (a), $e=0.1$, $\omega=1.8$ in panel (b), $e=0$,
   $\omega=0.87$ in panel (c), and $e=0.1$, $\omega=0.87$ in panel
   (d). Black dots are the attractors computed for the initial
   condition $\theta_0=0.01\pi$, $\dot\theta_0=0$ with Dynamics
   \cite{Dynamics}. The stable (thick dark green) and unstable (thin
   light red) periodic orbits, computed with (R)AUTO \cite{auto,S07},
   are superimposed.  Oscillations are shown in the $\theta$-range of
   $[-\pi,\pi]$, rotations are shown in the $\theta$-range
   $[\pi,3\pi]$. The damping $\gamma$ is $0.1$.  The basins of
   attraction shown in figure~\ref{fig:basins} have been computed at
   $p=0.5$, $p=1$, $p=1.4$ (see vertical dashed lines in panel (a) and
   (b)).  }
     \label{fig:oneparabif1}
\end{figure}

Figure~\ref{fig:oneparabif1} shows a series of four
one-parameter bifurcation diagrams for varying forcing amplitude $p$.
We pick two values for the frequency (the same as in
\cite{RefWorks:49}): $\omega=1.8$ (panel (a) and (b)), which is in the
1:2 resonance tongue, and $\omega=0.87$ (panel (c) and (d)), which is
in the 1:1 resonance tongue, and two values for the ellipticity: $e=0$
(panel (a) and (c)) and $e=0.1$ (panel (b) and (d)). The vertical axis
of all panels shows the coordinate $\theta$ of the stroboscopic map of
\eqref{eq:elliptic1a} taken at $t=2k\pi/\omega$ where $k$ is a large
integer.  Stable oscillations and rotations are dark green thick
lines, unstable oscillations and rotations are bright red thin
lines. All bifurcation curves in the two-parameter diagrams
figure~\ref{fig:paraspaces_oscs} and figure~\ref{fig:paraspaces_rots}
have been constructed by continuing the bifurcations shown as dark
circles in figure~\ref{fig:oneparabif1}. We have shifted the value of
$\theta$ by $2\pi$ for all rotations to prevent curves associated with
rotations and oscillations obscuring one another.  The underlying
black dots show the long-time behavior from the initial conditions
$(\theta_0,\dot\theta_0)=(0.01\pi,0)$ (and
$(\theta_0,\dot\theta_0)=(2\pi+0.01\pi,0)$) after waiting for a
transient of 1000 periods of excitation, computed with Dynamics
\cite{Dynamics}.

The main feature of the transition from $\epsilon=0$ to nonzero
$\epsilon$ is the perturbation of the reflection symmetry. The
symmetric system (with $\epsilon=0$) has the pitchfork bifurcations A
(for period two in Figure~\ref{fig:oneparabif1}(a)) and K$_1$ (for
period one in Figure~\ref{fig:oneparabif1}(c)) linking families of
symmetric and nonsymmetric oscillations. These pitchfork bifurcations
are perturbed into fold bifurcations (also named A and K$_1$ in the
figures \ref{fig:oneparabif1}(b) and (d)).

The rotations (which are nonsymmetric orbits) and the nonsymmetric
oscillations come in pairs of orbits symmetric to each other and lying
on top of each other for $\epsilon=0$ in figures
\ref{fig:oneparabif1}(a) and (c). The same applies to the bifurcations
of the nonsymmetric orbits: the fold J and the period doubling G of
the rotations, and the period doublings F and L for the nonsymmetric
oscillations (starting rapidly accumulating period doubling cascades)
are symmetric pairs of bifurcations, occuring simultaneously. This
symmetry is broken by the increase of $\epsilon$ such that the
formerly symmetric branches are now different: rotations in the
negative direction emerging from the fold J$^\mathrm{n}$ already exist
for smaller forcing $p$ than the rotations in the positive direction
emerging from J$^\mathrm{p}$, which are shifted toward larger forcing
$p$.  Similarly, the formerly symmetric pairs of nonsymmetric
oscillations lose their symmetry: one family is always shifted toward
larger $p$ (born at the fold F$_2$ in
Figure~\ref{fig:oneparabif1}\,(b), and K$_2$ in
Figure~\ref{fig:oneparabif1}\,(d)), the other family becomes a
continuous extension of the formerly symmetric oscillation.

For $\omega=1.8$ the visibility of chaotic attractors (bands of small
black dots in figure~\ref{fig:oneparabif1}) is shifted toward larger
$p$ by the symmetry breaking because stable periodic rotations exist
for larger $p$ (up to $p\approx 1.35$). At $\omega=0.87$ the
simulation also showed period-two oscillations jumping between two
potential wells for larger $p$ in the simulation results (black lines
evident after the chaotic bands in panels (c) and (d)).

\begin{figure}[ht]
 \begin{center}
   {\includegraphics[width=1.0\picrefwidth]{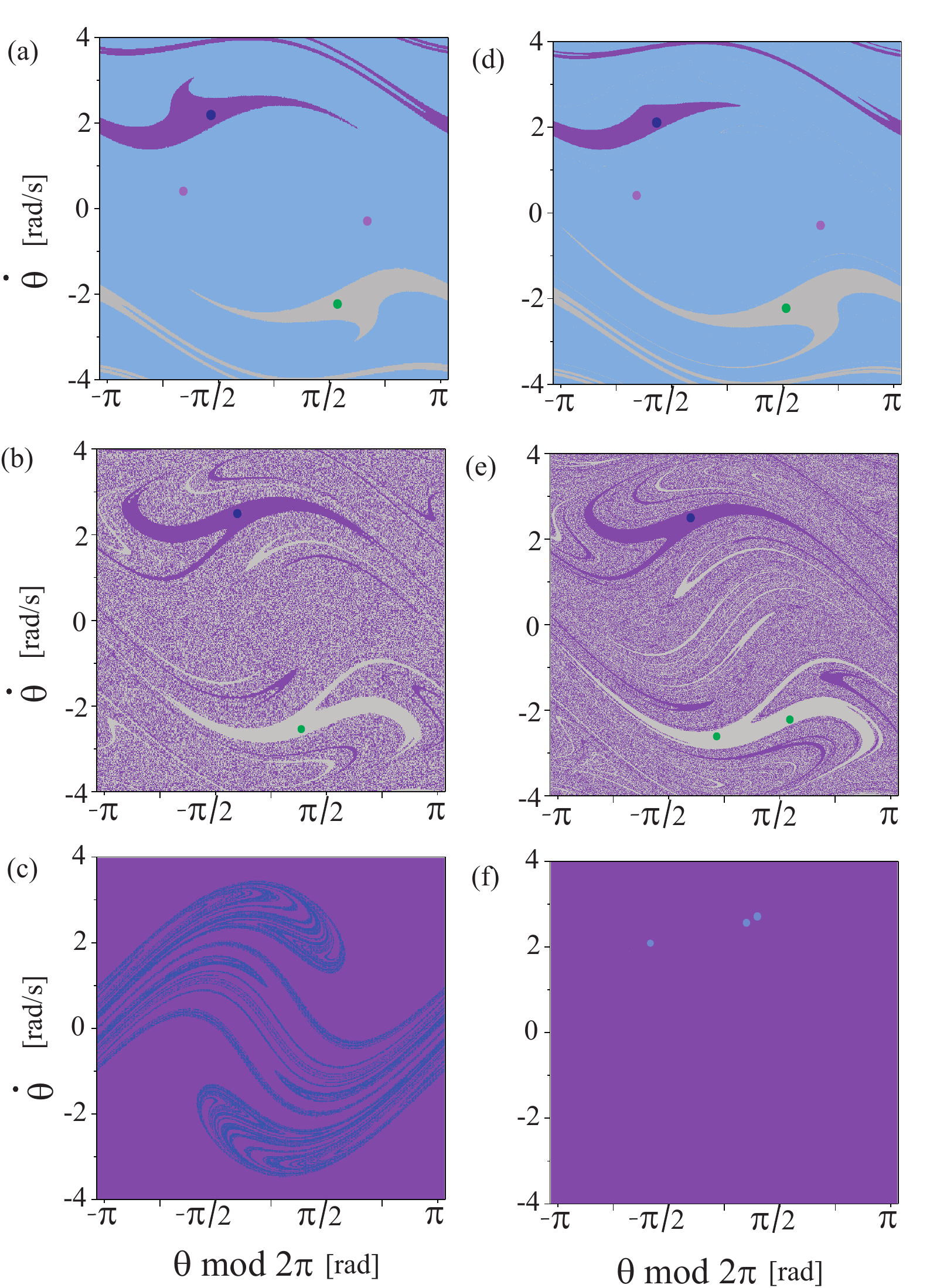}}
     \end{center}
     \caption{(colour online) Basins of attraction for different
       values of the forcing amplitude, $p$, and ellipticity
       $\epsilon$, computed with Dynamics \cite{Dynamics}. Panel (a):
       $\epsilon=0$, $p=0.5$, (b): $\epsilon=0$, $p=1$, (c):
       $\epsilon=0$, $p=1.4$, (d): $\epsilon=0.1$, $p=0.5$, (e):
       $\epsilon=0.1$, $p=1$, (f) $\epsilon=0$, $p=1.4$. The other
       parameters are $\gamma=0.1$ and $\omega =1.8$.  Stable periodic
       orbits are shown as large dots (their colour is chosen to give
       a contrast to their own basin of attraction).}
     \label{fig:basins}
\end{figure}

Figure~\ref{fig:basins} shows how the basins of attraction lose their
symmetry when one increases $\epsilon$ from $0$ to $0.1$.  The colour
coding of each point in the $(\theta,\dot\theta)$-plane is chosen
according to the attractor which the stroboscopic map reaches starting
from this point.
The forcing is $p=0.5$ in figure~\ref{fig:basins}\,(a) and (d), $p=1$
in panels (b) and (e), and $p=1.4$ in panels (c) and (f). The periodic
attractors (shown as dots) with $\dot\theta$-coordinates $\approx\pm2$
correspond to periodic rotations. The panels (d) and (e) show that for
small $p$ ($p=0.5$ and $p=1$) one direction of rotation (negative) has
a visibly larger basin of attraction than the other. At $p=1.4$ the
change of $\epsilon$ caused a crisis of the chaotic attractor in panel
(c), creating a period-three rotation.

For further increase of $\epsilon$ the effect that one attractor of the
formerly symmetric pair of nonsymmetric periodic orbits is shifted
toward higher values of forcing becomes more pronounced (as shown in
the the two-parameter diagrams in figure~\ref{fig:paraspaces_oscs}(c)
and figure~\ref{fig:paraspaces_rots}(b)). This shift depends
strongly and nonlinearly on $\epsilon$: for example, formula \eqref{eq:foldapp}
already underestimates this shift for J$^\mathrm{p}$ by 20\% for
$\epsilon=0.5$, which still corresponds to a narrow ellipse.

\section{Conclusions}
\label{sec:conclusion} Introducing a horizontal component into the
excitation of the classical parametrically excited pendulum results in
a symmetry breaking scenario. The excitation changes from a purely
vertical motion to a motion on an ellipse. The main effects of this
ellipticity are twofold: first, the well-known 1:2 and 1:1 resonance
tongues of the classical parametric pendulum merge into a single
region of instability, bounded by a period doubling and a fold
(saddle-node) bifurcation of the small amplitude oscillation.  Second,
rotations of the pendulum that have the same direction as the base
motion pick up energy from the horizontal excitation such that they
are present at lower overall forcing amplitudes. For example,
clockwise motion of the pivot around the ellipse results in a
preference for clockwise rotations of the pendulum.

Both effects of ellipticity are favorable for rotation:
the first effect implies that small-amplitude oscillations
around the hanging-down position, which are attractors competing with
rotations, lose their stability for smaller forcing
compared with the $e=0$ case. 
The second effect means that the parameter region
in the frequency-amplitude plane where rotations are supported
increases with increasing ellipticity. The excitation amplitude necessary to sustain
rotations in both directions also increases with increasing $\epsilon$ because of the
increasing `imperfection' of the symmetry.


A comparison between the bifurcation scenarios of the model and an
actual experiment is still outstanding. Direct bifurcation analysis
for experiments is a challenging task that may require the development
of entirely new experimental methods. Apart from this lack of
experimental verification, other open questions are: the
high-frequency approximation \eqref{eq:foldapp} suggests that for a
circular excitation ($\epsilon=1$, $\alpha=\pi/4$) rotations against
the base excitation are impossible regardless of the level of forcing
and damping. This is not true in general for lower frequency and
sufficiently small damping. Thus, we expect that, depending on the
shape of the excitation, there must be a critical damping level below
which rotations against the excitation direction become possible for a
suitable range of the excitation amplitudes and frequencies.

The small dissipation restricts the type of bifurcations and regimes
encountered in the system (for example, torus bifurcations are
impossible). We expect that even a small amount of interaction between
the pendulum and the base (\cite{xuzamm}) will lead to large regions
in the frequency-forcing plane where one can observe
quasi-periodicity.  Escape from a potential well tends to lead to
indeterminacy as introduced in \cite{RefWorks:61}. The precise
sequence of heteroclinic tangencies leading to escape from the
potential well is still largely unknown even for the parametrically
excited pendulum.

\section*{Acknowledgments}
\label{sec:ack}
  M.~W. acknowledges financial support by The Royal Society. B.~H. would like to thank EPSRC for financial support
  throughout his doctoral studies.

\bibliographystyle{elsarticle-num}

\end{document}